\documentclass[a4paper]{article}
\usepackage{amsmath,amssymb,amscd,amsthm}
\usepackage{epic,eepic}
\usepackage[dvips]{graphics}
\usepackage[usenames]{color}
\usepackage{comment}
\usepackage{enumerate}

\DeclareFontEncoding{OT2}{}{}
\DeclareFontSubstitution{OT2}{wncyr}{m}{n}
\DeclareSymbolFont{cyss}{OT2}{wncyss}{m}{n}
\DeclareSymbolFont{cyr}{OT2}{wncyr}{m}{n}
\DeclareMathSymbol{\sh}{\mathbin}{cyss}{`x}

\definecolor{gray}{gray}{0.5}

\newcommand{\End}{\operatorname{End}}
\newcommand{\ad}{\operatorname{ad}}
\newcommand{\PGL}{\operatorname{PGL}}

\newcommand{\Li}{\operatorname{Li}}

\newcommand{\id}{\operatorname{id}}

\newcommand{\C}{{\mathbf C}}
\newcommand{\R}{{\mathbf R}}
\newcommand{\Q}{{\mathbf Q}}

\newcommand{\bP}{{\mathbf P}}
\newcommand{\bunit}{{\mathbf I}}
\newcommand{\bnull}{{\mathbf 1}}

\newcommand{\fS}{{\mathfrak S}}
\newcommand{\fX}{{\mathfrak X}}

\newcommand{\fM}{{\mathfrak M}}
\newcommand{\fA}{{\mathfrak A}}

\newcommand{\cB}{{\mathcal B}}

\newcommand{\cL}{{\mathcal L}}
\newcommand{\cM}{{\mathcal M}}
\newcommand{\cU}{{\mathcal U}}
\newcommand{\cW}{{\mathcal W}}
\newcommand{\hcL}{{\hat{\mathcal L}}}
\newcommand{\tcL}{{\tilde{\mathcal L}}}

\theoremstyle{definition}
\newtheorem{thm}{Theorem}

\newtheorem{prop}[thm]{Proposition}

\newtheorem{rem}[thm]{Remark}

\title{\huge Iterated integrals and relations of multiple polylogarithms\\[0.5cm]}
\author{OI, Shu\thanks{Department of Mathematics, School of Fundamental Sciences and 
Engineering, Faculty of Science and Engineering, Waseda university. \endgraf \hspace{1em}3-4-1, Okubo, Shinjuku-ku, 
Tokyo 169-8555, Japan. \endgraf \hspace{1em}e-mail: {\tt shu\_oi@toki.waseda.jp}}
\and UENO, Kimio\thanks{Department of Mathematics, School of Fundamental Sciences 
and Engineering, Faculty of Science and Engineering, Waseda university. \endgraf \hspace{1em}3-4-1, Okubo, Shinjuku-ku, 
Tokyo 169-8555, Japan. \endgraf \hspace{1em}e-mail: {\tt uenoki@waseda.jp}
}}
\date{}

\begin{document}

\allowdisplaybreaks

\maketitle

\begin{abstract}
This is a summary for the authors' article \cite{OU1}
(prerint (2009) arXiv: math.QA/0910.0718), including a new result on the five term 
relation for the dilogarithm. This note will appear in the RIMS K\^{o}ky\^{u}roku for the 
conference on ``Representation Theory and Combinatorics'' held at Hokkaido University 
from August 25th to 28th, 2009. 
\end{abstract}

\section{Introduction}

The aim of our work is to construct and research the fundamental solution of the formal 
KZ (Knizhnik-Zamolodchikov) equation via iterated integrals. 
First we establish the decomposition theorem for the normalized fundamental solution of 
the formal KZ equation on the moduli space $\cM_{0,5}$ 
(or, the formal KZ equation of two variables).
Next we show that, 
by using iterated integrals, it can be viewed as a generating function of hyperlogarithms 
of the type $\cM_{0,5}$. The decomposition theorem says that the normalized fundamental 
solution decomposes to a product of two factors which are the normalized fundamental 
solutions of the formal (generalized) KZ equations of one variable. Comparing the different ways 
of decomposition gives the generalized harmonic product relations of the hyperlogarithms.  
These relations properly contain the harmonic product of multiple polylogarithms. 

The most simple case of the harmonic product is the following: Let us define 
\begin{align*}
\Li_{k_1,\ldots,k_r}(z)&=\sum_{n_1>\cdots>n_r>0}\frac{z^{n_1}}{n_1^{k_1}\cdots n_r^{k_r}},\\
\Li_{k_1,\ldots,k_{i+j}}(i,j;z_1,z_2)&=\sum_{n_1>\cdots>n_{i+j}>0}
                  \frac{z_1^{n_1} z_2^{n_{i+1}}}{n_1^{k_1}\cdots n_{i+j}^{k_{i+j}}}.
\end{align*}
Then we obtain 
\begin{align}
\Li_{k}(z_1)\Li_{l}(z_2)&=\sum_{m>0}\frac{z_1^m}{m^k} \sum_{n>0}\frac{z_2^n}{n^l}
=\left(\sum_{m>n>0}+\sum_{m=n>0}+\sum_{n>m>0}\right)\frac{z_1^m z_2^n}{m^k n^l} \nonumber \\
&=\Li_{k,l}(1,1;z_1,z_2)+\Li_{k+l}(z_1z_2)+\Li_{l,k}(1,1;z_2,z_1). \tag{HPMPL} \label{hpmpl}
\end{align}
Taking the limit, we have the harmonic product of multiple zeta values 
\begin{align}\tag{HPMZV} \label{hpmzv}
   \zeta(k) \zeta(l) = \zeta(k,l)+\zeta(k+l)+\zeta(l,k). 
\end{align}
(The harmonic product of multiple zeta values is considered from the viewpoint of 
arithmetic geometry in \cite{BF},\  \cite{DT},\  \cite{F}.)

Moreover we consider the transformation theory of the fundamental solution
of the formal KZ equation of two variables 
and derive the five term relation for the dilogarithm due to Hill \cite{Le}, 
\begin{align}
\Li_2(z_1z_2)= \Li_2 & \left(\frac{-z_1(1-z_2)}{1-z_1}\right)  + 
\Li_2\left(\frac{-z_2(1-z_1)}{1-z_2}\right)  \nonumber \\
& + \Li_2(z_1)+\Li_2(z_2)+\frac{1}{2}\log^2\left(\frac{1-z_1}{1-z_2}\right). 
\tag{5TERM} \label{5term}
\end{align}

For detailed accounts of the results in this note, see \cite{OU1} and \cite{OU2}. 
The transformation theory of the formal KZ equation of one variable 
(or the formal KZ equation on $\cM_{0,4}$) is studied in \cite{OkU}.

\paragraph{Acknowledgment}

The authors express their gratitude to Professor Hideaki Morita for giving them 
a chance of a lecture. 
The second author is partially supported by JPSP Grant-in-Aid No. 19540056.

\section{The formal KZ equation on $\cM_{0,n}$}
\subsection{Definition of the formal KZ equation}

First we introduce the formal KZ equation: It is defined on the configuration space 
of $n$ points of $\bP^1$ ($=$ the complement of the hyperplane arrangement associated 
with Dynkin diagram of $A_{n-1}$-type), which is by definition 
\begin{align*}
(\bP^1)^n_*=\{(x_1,\ldots,x_n) \in 
\underbrace{\bP^1\times\cdots\times\bP^1}_{n} \;|\; x_i\neq x_j\;\; (i\neq j)\}.
\end{align*}
The \textbf{infinitesimal pure braid Lie algebra} 
\begin{align*}
\fX=\fX(\{X_{ij}\}_{1\le i,j \le n}) := \C\{X_{ij}\;|\;{1\le i,j \le n}\}
                                    \Big/ \eqref{ipbr}
\end{align*}
is a graded Lie algebra for the lower central series of the fundamental group 
of $(\bP^1)^n_*$ \ \cite{I}. It is generated by the formal elements 
$\{X_{ij}\}_{1\le i,j \le n}$ with the defining relations \eqref{ipbr} 
(the infinitesimal pure braid relations) 
\begin{align}
\begin{cases}
X_{ij}=X_{ji}, \quad &X_{ii}=0,\\
\sum_j X_{ij}=0 \quad (\forall i ),  
\quad &[X_{ij},X_{kl}]=0 \quad (\{i,j\} \cap \{k,l\} = \emptyset). 
\tag{IPBR}\label{ipbr}
\end{cases}
\end{align}
By $\cU(\fX)$, we denote the universal enveloping algebra of $\fX$. It has the unit $\bunit$ 
and has the grading with respect to the homogeneous degree of an element: 
\begin{align*}
          \cU(\fX)= \bigoplus_{s=0}^{\infty} \, \cU_s(\fX).
\end{align*}

The formal KZ equation is by definition 
\begin{align}\tag{KZ}\label{KZeq}
dG=\varOmega G,  \qquad
   \varOmega = \sum_{i<j}\xi_{ij}X_{ij}, \qquad \xi_{ij}=d\log(x_i-x_j), 
\end{align}
which is a $\fX$-valued total differential equation (or, a connection) on $(\bP^1)^n_*$. 
(Such a formal equation was considered in \cite{Ha},\ \cite{De},\ \cite{Dr},\ \cite{W}.)

The 1-forms $\xi_{ij}$'s satisfy only the \textbf{Arnold relations} \cite{A} 
as non-trivial relations of degree 2: 
\begin{align}\tag{AR}\label{AR}
      \xi_{ij} \wedge \xi_{ik} + \xi_{ik} \wedge \xi_{jk} + 
      \xi_{jk} \wedge \xi_{ij} = 0.
\end{align}
From \eqref{ipbr} and \eqref{AR}, one can see that \eqref{KZeq} is integrable and 
has $\mathrm{PGL}(2,\C)$-invariance. 
Hence \eqref{KZeq} can be viewed as an equation on the moduli space 
\begin{align*}
\cM_{0,n}=\PGL(2,\C) \backslash (\bP^1)^n_*.
\end{align*}
Hereafter we will call \eqref{KZeq} the \textbf{formal KZ equation on the moduli space} $\cM_{0,n}$.

\subsection{The formal KZ equation on $\cM_{0,4}$ and $\cM_{0,5}$}

For analysis of \eqref{KZeq}, it is convenient to use the \textbf{cubic coordinates} on 
$\cM_{0,n}$ \cite{B}. Introducing the simplicial coordinates $\{y_i\}$ by 
\begin{align*}
y_i=\frac{x_i-x_{n-2}}{x_i-x_n}\frac{x_{n-1}-x_{n}}{x_{n-1}-x_{n-2}}
                                                \qquad (i=1,\ldots,n-3), 
\end{align*}
(fixing three points $y_n=\infty,\; y_{n-1}=1,\; y_{n-2}=0$) 
the cubic coordinates $\{z_i\}$ are defined by blowing up at the origin, 
\begin{align*}
   y_i \,=\, z_1 \cdots z_i  \qquad (i=1,\dots,n-3).
\end{align*}

We give representations of \eqref{KZeq} for $n=4, 5$. In the cubic coordinates of $\cM_{0,4}$, 
we put $z=z_1$ and $Z_1=X_{12}, Z_{11}=-X_{13}$. Then \eqref{KZeq} is represented as 
\begin{align}\tag{1KZ}\label{1KZeq}
dG=\varOmega G, \qquad 
\varOmega=\zeta_1Z_1+\zeta_{11}Z_{11}, \qquad 
\zeta_1=\frac{dz}{z},\; \zeta_{11}=\frac{dz}{1-z}. 
\end{align}
which is referred to as the \textbf{formal KZ equation of one variable}. 
The singular divisors of this equation are $D(\cM_{0,4}^{cubic}):=\{z=0,1,\infty\}$. 
The Lie algebra $\fX$ is a free Lie algebra generated by $Z_1, Z_{11}$, and 
\eqref{AR} reduces to the trivial one $\zeta_1 \wedge \zeta_{11}=0$.

In the case of $\cM_{0,5}$, we put 
\begin{align*}
Z_1=X_{12}+X_{13}+X_{23},\; Z_{11}=-X_{14},\; Z_2=X_{23},\; Z_{22}=-X_{12},\ Z_{12}=-X_{24}. 
\end{align*}
In the cubic coordinates of $\cM_{0,5}$, (KZ) reads as 
\begin{gather}\tag{2KZ}\label{2KZeq}
dG = \varOmega G, \qquad 
\varOmega=\zeta_1 Z_1 + \zeta_{11} Z_{11} + \zeta_2 Z_2 + \zeta_{22} Z_{22} + \zeta_{12} Z_{12}, \\
\zeta_1=\frac{dz_1}{z_1},\;\; \zeta_{11}=\frac{dz_1}{1-z_1},\;\; 
\zeta_2=\frac{dz_2}{z_2},\;\; \zeta_{22}=\frac{dz_2}{1-z_2},\;\; 
\zeta_{12}=\frac{d(z_1z_2)}{1-z_1z_2}, \notag
\end{gather}
which is referred to as the \textbf{formal KZ equation of two variables}. 
The singular divisors of this equation are 
$D(\cM_{0,5}^{cubic}):=\{z_1=0,1,\infty\}\cup\{z_2=0,1,\infty\}\cup\{z_1z_2=1\}$. 
The Lie algebra $\fX$ is generated by the five elements $Z_1, Z_{11}, Z_2, Z_{22}, Z_{12}$ 
with the defining relations 
\begin{align}
\begin{cases}
[Z_1,Z_2]=[Z_{11},Z_2]=[Z_1,Z_{22}]=0, \\
[Z_{11},Z_{22}]=[-Z_{11},Z_{12}]=[Z_{22},Z_{12}]=[-Z_1+Z_2,Z_{12}].
\end{cases} \tag{IPBR'}\label{ipbr'}
\end{align}
Non trivial relations among \eqref{AR} are 
\begin{align}
\begin{cases}
(\zeta_1+\zeta_2) \wedge \zeta_{12} = 0,\\
\zeta_{11} \wedge \zeta_{12} + \zeta_{22} \wedge (\zeta_{11} - \zeta_{12}) 
                      - \zeta_2 \wedge \zeta_{12}=0.
\end{cases} \tag{AR'}\label{AR'}
\end{align}

The following is a figure of the divisors $D(\cM^{cubic}_{0,5})$. 
Note that they are normal crossing at $(z_1,z_2)=(0,0),(1,0),(0,1)$.
\begin{center}
\setlength{\unitlength}{1cm}
\begin{picture}(6,5.5)(0,0)
\put(0,0){\scalebox{0.5}{\includegraphics{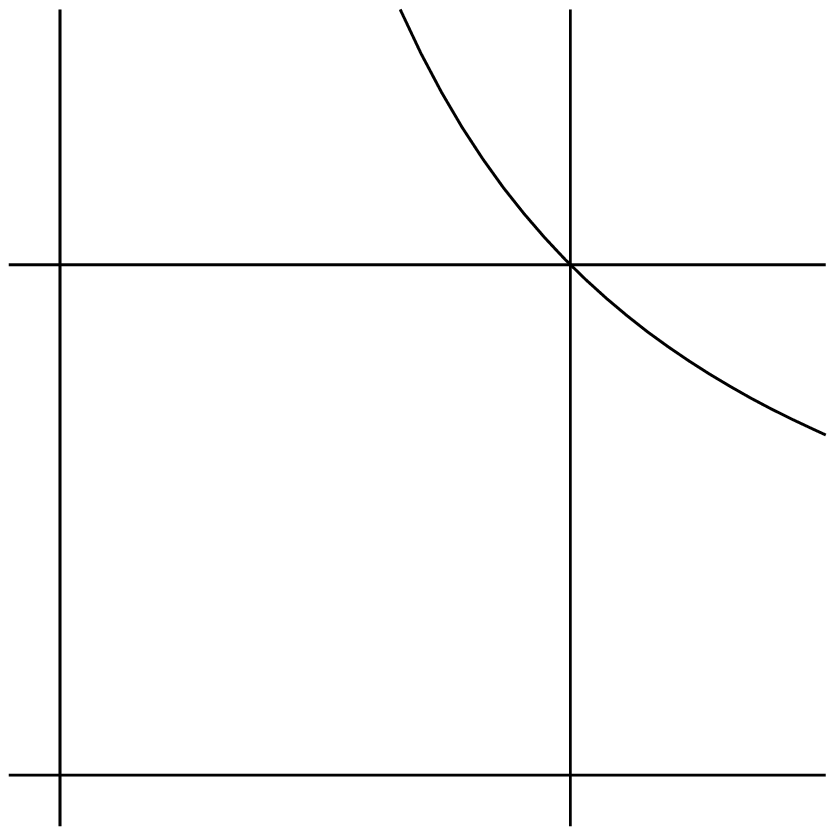}}}
\put(-0.5,0){$(0,0)$}
\put(3.1,0){$(1,0)$}
\put(-0.5,3.1){$(0,1)$}
\put(3.1,3.1){$(1,1)$}
\put(4.4,0.3){$z_1$}
\put(0.3,4.4){$z_2$}
\end{picture}
\end{center}

\section{The fundamental solution of the formal KZ equation on $\cM_{0,4}$}
\subsection{A free shuffle algebra and iterated integral on $\cM_{0,4}$}

For a free shuffle algebra $S=S(a_1,\ldots,a_r)$ generated by the alphabet $a_1,\ldots,a_r$, 
we denote by $\bnull$ the unit, by $\circ$ the product of concatenation and by $\sh$ 
the shuffle product:
\begin{gather*}
S=(\C\langle a_1,\ldots,a_r\rangle, \sh),\\
w \sh \bnull = \bnull \sh w = \bnull,\\
(a_i \circ w) \sh (a_j \circ w')=a_i \circ (w \sh (a_j \circ w')) 
+ a_j \circ ((a_i \circ w) \sh w').
\end{gather*}
It is a graded algebra with respect to the homogeneous degree of an element.

Let $\zeta_1, \zeta_{11}$ be the 1-forms in \eqref{1KZeq}, and $S(\zeta_1,\zeta_{11})$ 
a free shuffle algebra generated by them. For any word 
$\varphi=\omega_1 \circ \cdots \circ \omega_r \quad (\omega_i=\zeta_1, \mbox{or}, \zeta_{11})$ 
in $S(\zeta_1,\zeta_{11})$, we set the iterated integral by 
\begin{align*}
\int_{z_0}^z \,\, \varphi
=\int_{z_0}^z \,\, \omega_1(z') \, \int_{z_0}^{z'}\,\,\omega_2 \circ \cdots \circ \omega_r, 
\end{align*}
which gives a many-valued analytic function on $\bP^1-D(\cM^{cubic}_{0,4})$.

For $\varphi,\psi \in S(\zeta_1,\zeta_{11})$, we have 
\begin{align*}
\int (\varphi \sh \psi)=\left(\int \varphi \right)\left(\int \psi \right). 
\end{align*}

A free shuffle algebra has the structure of a Hopf algebra, and  
$S(\zeta_1,\zeta_{11})$ is a dual Hopf algebra of the universal enveloping 
algebra $\cU(\fX)$.

\subsection{The fundamental solution of \eqref{1KZeq}}

Next we consider the fundamental solution of \eqref{1KZeq} normalized at the origin $z=0$. 
We denote it by  $\cL(z)$. It is a solution satisfying the following condition: 
\begin{align*}
\cL(z)=\hcL(z)z^{Z_1}
\end{align*}
where $\hcL(z)$ is represented as 
\begin{align*}
\hcL(z)= \sum_{s=0}^{\infty} \hcL_s(z), \quad  \hcL_s(z) \in \cU_s(\fX), 
                       \quad \hcL_s(0)=0\;\; (s>0), \quad \hcL_0(z)=\bunit.
\end{align*}
It is easy to see that $\hcL_s(z)$ satisfies the following recursive equation:
\begin{align*}
\frac{d\hcL_{s+1}}{dz} = \frac{1}{z}[Z_1, \hcL_{s}] + \frac{1}{1-z}Z_{11} \hcL_{s} 
\quad (s=0,1,2,\dots).
\end{align*}
Since the term $\frac{1}{z}[Z_1, \hcL_{s}]$ is holomorphic at $z=0$, $\hcL_{s+1}(z)$ 
is uniquely determined by 
\begin{align*}
\hcL_{s+1}(z)= \int_0^z \Big( \frac{1}{z}[Z_1, \hcL_{s}] + \frac{1}{1-z}Z_{11}\hcL_{s} \Big) \, dz.
\end{align*}
In terms of iterated integral, it is expressed as 
\begin{align*}
\hcL_s(z) &=\sum_{k_1+\cdots+k_r=s} 
            \Big\{\int_0^z \zeta_1^{k_1-1}\circ \zeta_{11} \circ\cdots\circ \zeta_1^{k_r-1} 
                                                         \circ \zeta_{11}\Big\} \notag \\
& \hphantom{\sum_{k_1+\cdots+k_r=s} \Big\{} \times \ad(Z_1)^{k_1-1}\mu(Z_{11})
                                   \cdots\ad(Z_1)^{k_r-1}\mu(Z_{11})(\bunit).
\end{align*}
Here $\ad(Z_1) \in \End(\cU(\fX))$ stands for the adjoint operator by $Z_1$, 
and $\mu(Z_{11}) \in \End(\cU(\fX))$ the multiplication of $Z_{11}$ from the left.
From these considerations, it follows that \textbf{the fundamental solution normalized 
at $z=0$ exists and is unique}. 

The iterated integral in the right hand side is a \textbf{multiple polylogarithm 
of one variable}:
\begin{align}\tag{1MPL} \label{1mpl}
\Li_{k_1,\dots,k_r}(z) = \int_0^z \zeta_1^{k_1-1}\circ \zeta_{11} 
                                          \circ\cdots\circ \zeta_1^{k_r-1} \circ \zeta_{11}. 
\end{align}
If $|z|<1$, it has a Taylor expansion
\begin{align*}
\Li_{k_1,\dots,k_r}(z)=\sum_{n_1>n_2>\cdots>n_r>0} \frac{z^{n_1}}{n_1^{k_1}\cdots n_r^{k_r}}. 
\end{align*}
If $k_1 \geq 2$, we have 
\begin{align*}
\lim_{z \to 1-0}\, \Li_{k_1,\dots,k_r}(z) = \zeta(k_1,\dots,k_r),
\end{align*}
where the right side above is a \textbf{multiple zeta value}, 
\begin{align}\tag{MZV}\label{mzv}
\zeta(k_1,\dots,k_r)=\sum_{n_1>\cdots>n_r>0} \frac{1}{n_1^{k_1}\cdots n_r^{k_r}}.
\end{align}

\subsection{The fundamental solution of the formal generalized KZ equation of one variable}

Let us consider a generalization of \eqref{1KZeq}. 
For mutually distinct points $a_1, \ldots, a_m \in \C-\{0\}$ we set
\begin{equation}\tag{G1KZ} \label{g1KZeq}
dG= \varOmega G, \quad \varOmega=\frac{dz}{z}X_0+\sum_{i=1}^{m}\frac{a_i dz}{1-a_i z}X_i.
\end{equation}
Here the coefficients $X_0, X_1,\ldots, X_m$ are free formal elements.
For $r=1, a_1=1$, this is the formal KZ equation of one variable. 
This is a differential equation of the Schlesinger type 
with regular singular points $0,1/a_1,\ldots, 1/a_m, \infty$.
We call \eqref{g1KZeq} the \textbf{formal generalized KZ equation of one variable}. 

Let $\fX=\C\{X_0,X_1,\ldots,X_m \}$ be a free Lie algebra generated by 
$X_0,X_1,\dots,X_m$, and $\cU(\fX)$  the universal enveloping algebra.

The free shuffle algebra
$S(\xi_0, \xi_1,\ldots, \xi_m)$ where 
\begin{align*}
\xi_0=\frac{dz}{z}, \quad\quad  \xi_i=\frac{a_idz}{1-a_iz}, \quad (1 \leq i \leq m),
\end{align*}
is a dual Hopf algebra of $\cU(\fX)$. 

\textbf{The fundamental solution $\cL(z)$ normalized at the origin $z=0$ of this equation 
exists and is unique}. It satisfies the following conditions:
\begin{align*}
  \cL(z)=\hcL(z)z^{X_0}
\end{align*}
where $\hcL(z)$ is represented as 
\begin{align*}
   \hcL(z)= & \sum_{s=0}^{\infty} \hcL_s(z), \quad  \hcL_s(z) \in \cU_s(\fX), 
                       \quad \hcL_s(0)=0\;\; (s>0), \quad \hcL_0(z)=\bunit, \\*
\hcL_s(z) &= \sum_{\substack{k_1+\cdots+k_r=s\\i_1,\ldots,i_r \in \{1,\ldots,m\}}} 
L({}^{k_1}a_{i_1}\cdots{}^{k_r}a_{i_r};z) \notag \\*
& \hphantom{\sum_{i_1,\ldots,i_r \in \{1,\ldots,m\}}} \times 
  \ad(X_0)^{k_1-1}\mu(X_{i_1})\cdots\ad(X_0)^{k_r-1}\mu(X_{i_r})(\bunit). 
\end{align*}
Here $L({}^{k_1}a_{i_1}\cdots{}^{k_r}a_{i_r};z)$ is a 
\textbf{hyperlogarithm of the general type}: 
\begin{align}\tag{HLOG} \label{hlog}
             L({}^{k_1}a_{i_1}\cdots{}^{k_r}a_{i_r};z) 
               := \int_0^z \xi_0^{k_1-1}\circ \xi_{i_1} \circ \xi_0^{k_2-1}\circ 
                               \xi_{i_2} \circ\cdots\circ \xi_0^{k_r-1} \circ \xi_{i_r}.
\end{align}
For $r=1$ and $a_1=1$, this is \eqref{1mpl}. 
If $|z|<\min\{\frac{1}{|a_{i_1}|},\ldots,\frac{1}{|a_{i_r}|}\}$, it has a Taylor expansion
\begin{align*}
L({}^{k_1}a_{i_1}\cdots{}^{k_r}a_{i_r};z)=
      \sum_{n_1>n_2>\cdots>n_r>0}\frac{a_{i_1}^{n_1-n_2}a_{i_2}^{n_2-n_3}\cdots 
                                     a_{i_r}^{n_r}}{n_1^{k_1}\cdots n_r^{k_r}} z^{n_1}.
\end{align*}

\section{The fundamental solution of the formal KZ equation on $\cM_{0,5}$}
\subsection{The reduced bar algebra and iterated integrals on $\cM_{0,5}$}

Let $S=S(\zeta_1,\zeta_{11},\zeta_2,\zeta_{22},\zeta_{12})$ be a free shuffle algebra 
generated by $\zeta_1,\zeta_{11},\zeta_2,\zeta_{22},\zeta_{12}$ which are 1-forms in \eqref{2KZeq}. 
The iterated integral of an element in S, in general, depends on the integral path. 
We want to construct a shuffle subalgebra of S such that the iterated integral of any 
element in this subalgebra depends only on the homotopy class of the integral path. 
We say that an element
\begin{align*}
S \ni \varphi \,=\, \sum_{I=(i_1,\ldots,i_s)} \ c_I \omega_{i_1}\circ\cdots\circ\omega_{i_s}, 
\end{align*}
where $\omega_{i} \in \{\zeta_1,\zeta_{11},\zeta_2,\zeta_{22},\zeta_{12}\}$, satisfies 
\textbf{Chen's integrability condition} \cite{C1} if and only if 
\begin{align}\tag{CIC}\label{cic}
\sum_I c_I \; \omega_{i_1}\otimes\cdots\otimes\omega_{i_l}\wedge
              \omega_{i_{l+1}}\otimes\cdots\otimes\omega_{i_s}=0 
\end{align}
holds for any $l$ ($1 \le l<s$) as a multiple differential form.
Let $\cB$ be the subalgebra of elements satisfying \eqref{cic}. We call it 
the \textbf{reduced bar algebra}, which coincides with the 0-th cohomology of 
the reduced bar complex \cite{C2} associated with the Orlik-Solomon algebra \cite{OT} 
generated by $\zeta_1,\zeta_{11},\zeta_2,\zeta_{22},\zeta_{12}$. 

For any element $\varphi \in \cB$, the iterated integral 
\begin{align*}
    \int_{(z^{(0)}_{1},z^{(0)}_{2})}^{(z_1,z_2)} \ \varphi
\end{align*}
gives a many-valued analytic function on $\bP^1 \times \bP^1 - D(\cM_{0,5}^{cubic})$.

Let us consider more on the structure of $\cB$: It is a graded algebra; 
$\cB=\bigoplus_{s=0}^\infty \cB_s, \;\; \cB_s=\cB \cap S_s$ where $S_s$ denotes 
the degree s part of $S$: We have 
\begin{align*}
& \hspace{10mm} \cB_0=\C\bnull,\quad 
  \cB_1=\C\zeta_1\oplus\C\zeta_{11}\oplus\C\zeta_2\oplus\C\zeta_{22}\oplus\C\zeta_{12}, \\
& \hspace{10mm} \cB_2=\bigoplus_{\omega \in A}\C \omega\circ\omega \oplus 
\bigoplus_{i=1,2}\C \zeta_i\circ\zeta_{ii} \oplus \bigoplus_{i=1,2}\C \zeta_{ii}\circ\zeta_i\\*
& \hspace{15mm} \oplus \bigoplus_{\substack{\omega_1=\zeta_1,\zeta_{11}\\ \omega_2=\zeta_2,\zeta_{22}}}
\C(\omega_1\circ\omega_2+\omega_2\circ\omega_1) \oplus \bigoplus_{\omega \in A-\{\zeta_{12}\}}\C 
(\omega\circ\zeta_{12}+\zeta_{12}\circ\omega)\\*
& \hspace{15mm} \oplus \C (\zeta_1\circ\zeta_{12}+\zeta_2\circ\zeta_{12}) \oplus 
   \C (\zeta_{11}\circ\zeta_{12}+\zeta_{22}\circ\zeta_{11}-\zeta_{22}\circ\zeta_{12}
                                                                  -\zeta_2\circ\zeta_{12})
\end{align*}
where $A:=\{\zeta_1,\zeta_{11},\zeta_2,\zeta_{22},\zeta_{12}\}$.
For $s>2$, $\cB_s$ is characterized as follows \cite{B};
\begin{align*}
\cB_s&=\bigcap_{j=1}^{s-1}\cB_j\circ\cB_{s-j}
=\bigcap_{j=0}^{s-2}\underbrace{\cB_1\circ\cdots\circ\cB_1}_{j\text{ times}}
\circ\cB_2\circ\underbrace{\cB_1\circ\cdots\circ\cB_1}_{s-j-2\text{ times}}. 
\end{align*}

Put 
\begin{align*}
   \zeta_{12}^{(1)}= \frac{z_2dz_1}{1-z_1z_2}, \qquad 
                \zeta_{12}^{(2)}= \frac{z_1dz_2}{1-z_1z_2}.
\end{align*}
One can define a linear map 
\begin{align*}
   \iota_{1 \otimes 2 }:\cB \longrightarrow 
               S(\zeta_1,\zeta_{11},\zeta_{12}^{(1)}) \otimes S(\zeta_2,\zeta_{22})
\end{align*}
by the following procedure;
\begin{enumerate}
\item  \ pick up the terms only having a form $\psi_1 \circ \psi_2 \in 
      S(\zeta_1,\zeta_{11},\zeta_{12}) \circ S(\zeta_2,\zeta_{22})$.
\item \ change each term $\psi_1 \circ \psi_2$ to $\psi_1 \otimes \psi_2 \in 
    S(\zeta_1,\zeta_{11},\zeta_{12}) \otimes S(\zeta_2,\zeta_{22})$.
\item \ replace $\zeta_{12}$ to $\zeta_{12}^{(1)}$.
\end{enumerate}
A linear map 
\begin{align*}
 \iota_{2 \otimes 1 }:\cB \longrightarrow 
               S(\zeta_2,\zeta_{22},\zeta_{12}^{(2)}) \otimes S(\zeta_1,\zeta_{11})
\end{align*}
is defined in the same way. 

One can show that
\begin{align*}
\cU(\fX) & \;\cong\; \cU(\C\{Z_1,Z_{11},Z_{12}\})\otimes \cU(\C\{Z_2,Z_{22}\}) \\
         & \;\cong\; \cU(\C\{Z_2,Z_{22},Z_{12}\})\otimes \cU(\C\{Z_1,Z_{11}\}) 
\end{align*}
and that $\cB$ is a dual Hopf algebra of $\cU(\fX)$.
Through this isomorphism and the duality, one can show the following proposition:

\begin{prop}\label{prop:b}
The maps $\iota_{1 \otimes 2}$ and $\iota_{2 \otimes 1}$ are $\sh$-isomorphisms.
\end{prop}

(Such an isomorphism is also obtained by \cite{B}.) 

Let $\cB^0$ be the subspace of $\cB$ spanned by elements 
which have no terms ending with  $\zeta_1$ and $\zeta_2$, and 
$S^0(\zeta_1,\zeta_{11},\zeta_{12}^{(1)})$ (resp. $S^0(\zeta_2,\zeta_{22})$)
the subspace spanned by elements which have no terms ending with $\zeta_1$ 
(resp. $\zeta_2$), and so on. They are shuffle algebras. One can show the following 
isomorphism:

\begin{prop}\label{prop:b0}
By $\iota_{1 \otimes 2}$ and $\iota_{2 \otimes 1}$, 
\begin{align*}
    \cB^0 \, \cong \, S^0(\zeta_1,\zeta_{11},\zeta_{12}^{(1)})\otimes S^0(\zeta_2,\zeta_{22}) 
          \, \cong \, S^0(\zeta_2,\zeta_{22},\zeta_{12}^{(2)})\otimes S^0(\zeta_1,\zeta_{11}).
\end{align*}
\end{prop}

The free shuffle algebra $S(\zeta_1,\zeta_{11},\zeta_{12}^{(1)})$ is a polynomial 
algebra over $S^0(\zeta_1,\zeta_{11},\zeta_{12}^{(1)})$ of the variable $\zeta_1$ as a 
shuffle algebra \cite{R}:
\begin{align*}
        S(\zeta_1,\zeta_{11},\zeta_{12}^{(1)}) \, \cong \, 
        S^0(\zeta_1,\zeta_{11},\zeta_{12}^{(1)})[\zeta_1]. 
\end{align*}
Likewise, we have 
\begin{align*}
        S(\zeta_2,\zeta_{22}) \, \cong \, S^0(\zeta_2,\zeta_{22})[\zeta_2]
\end{align*}
as a shuffle algebra. Applying these isomorphisms to Proposition \ref{prop:b0}, we
have 

\begin{prop}\label{prop:structure_b}
The reduced bar algebra $\cB$ is a polynomial algebra over $\cB^0$ of the variables 
$\zeta_1, \zeta_2$ as a shuffle algebra: 
\begin{align*}
      \cB \, \cong \, \cB^0 [\zeta_1, \zeta_2].
\end{align*}
\end{prop}

Assume that $0<|z_1|,|z_2|<1$ and 
define the following two contours $C_{1 \otimes 2},\; C_{2 \otimes 1}$:

\setlength{\unitlength}{1.5cm}
\begin{picture}(0,5.4)(-0.2,-0.6)
\put(0,0){\scalebox{0.75}{\includegraphics{m05.eps}}}
\put(-0.5,0){$(0,0)$}
\put(3.1,0){$(1,0)$}
\put(-0.5,3.1){$(0,1)$}
\put(3.1,3.1){$(1,1)$}
\put(4.4,0.3){$z_1$}
\put(0.3,4.4){$z_2$}
\put(2.2,2.65){$(z_1,z_2)$}
\put(1.3,0.6){$C_{2\otimes1}^{(1)}$}
\put(1.9,1.3){$C_{2\otimes1}^{(2)}$}
\put(0.6,1.3){$C_{1\otimes2}^{(2)}$}
\put(1.3,2.2){$C_{1\otimes2}^{(1)}$}
\put(2.5,0.8){$C_{2\otimes 1}=C_{2\otimes1}^{(2)}\circ C_{2\otimes1}^{(1)}$}
\put(-0.4,2.6){$C_{1\otimes 2}=C_{1\otimes2}^{(1)}\circ C_{1\otimes2}^{(2)}$}

\thicklines
\put(0.55,0.5){\vector(1,0){2}}
\put(2.55,0.5){\vector(0,1){2}}
\thinlines
\put(0.55,0.5){\vector(0,1){2}}
\put(0.55,2.5){\vector(1,0){2}}
\put(4.5,2.8){\begin{minipage}{6.8cm}
$C_{1\otimes 2}=C_{1\otimes2}^{(1)}\circ C_{1\otimes2}^{(2)}$,\\ $\qquad C_{1\otimes2}^{(2)}\!:(0,0)
\to(0,z_2)$,\\ $\qquad C_{1\otimes2}^{(1)}:(0,z_2)\to(z_1,z_2)$.\\
$C_{2\otimes 1}=C_{2\otimes1}^{(2)}\circ C_{2\otimes1}^{(1)}$,\\ 
$\qquad C_{2\otimes1}^{(1)}:(0,0)\to(z_1,0)$,\\ $
\qquad C_{2\otimes1}^{(2)}:(z_1,0)\to(z_1,z_2) $.\\
\end{minipage}}
\end{picture}
\setlength{\unitlength}{1cm}

\noindent
The composition of paths $C \circ C'$ is defined by connecting $C$ after $C'$.

For $\psi_1 \otimes \psi_2 \in 
S^0(\zeta_1,\zeta_{11},\zeta_{12}^{(1)})\otimes S^0(\zeta_2,\zeta_{22}) $, we set 
\begin{align*}
\int_{C_{1\otimes 2}}\psi_1 \otimes \psi_2 := \int_{z_1=0}^{z_1}\psi_1 \int_{z_2=0}^{z_2}\psi_2
\end{align*}
and for $\psi_1 \otimes \psi_2 \in 
S^0(\zeta_2,\zeta_{22},\zeta_{12}^{(2)})\otimes S^0(\zeta_1,\zeta_{11})$, 
\begin{align*}
\int_{C_{2 \otimes 1}}\psi_1 \otimes \psi_2 := \int_{z_2=0}^{z_2}\psi_1 \int_{z_1=0}^{z_1}\psi_2.
\end{align*}
Since the map $\iota_{1 \otimes 2}$ (resp. $\iota_{2 \otimes 1}$) picks up the terms of $\cB^0$ 
whose iterated integral along $C_{1 \otimes 2}$ (resp. $C_{2 \otimes 1}$) does not vanish, 
we have
\begin{align*}
\int_{(0,0)}^{(z_1,z_2)} \!\!\varphi&=\int_{C_{1 \otimes 2}} \!\! \varphi
=\int_{C_{1 \otimes 2}} \!\!\iota_{1 \otimes 2}(\varphi) \\
&=\int_{C_{2 \otimes 1}} \!\! \varphi
=\int_{C_{2 \otimes 1}} \!\!\iota_{2 \otimes 1}(\varphi) 
\end{align*}
for $\varphi \in \cB^0$.

\subsection{The fundamental solution of \eqref{2KZeq}}

We consider the fundamental solution $\cL(z_1,z_2)$ of \eqref{2KZeq} normalized at 
the origin $(z_1,z_2)=(0,0)$. It is a solution satisfying the following conditions:
\begin{align*}
\cL(z_1,z_2)=\hcL(z_1,z_2) z_1^{Z_1}z_2^{Z_2}
\end{align*}
where
\begin{align*}
  \hcL(z_1,z_2)=\sum_{s=0}^{\infty}\hcL_s(z_1,z_2), \quad \hcL_s(z_1,z_2)\in\cU_s(\fX), \quad
 \hcL_s(0,0)=0 \ (s > 0),
\end{align*}
and $\hcL_0(z_1,z_2)=\bunit$. We put  
\begin{align*}
 \varOmega_0 &=\zeta_1 Z_1+\zeta_2 Z_2, \\
 \varOmega'  &=\varOmega-\varOmega_0=\zeta_{11} Z_{11}+\zeta_{22} Z_{22}+\zeta_{12} Z_{12}. 
\end{align*}
It is easy to see that $\hcL_s(z_1,z_2)$ satisfies the following recursive equation:
\begin{align*}
     d\hcL_{s+1}(z_1,z_2)=[\varOmega_0,\hcL_s(z_1,z_2)]+\varOmega'\hcL_s(z_1,z_2).
\end{align*}
Hence we have 
\begin{align}\tag{IISOL}  \label{iisol}
\hcL_s(z_1,z_2)=\int_{(0,0)}^{(z_1,z_2)} \left(\ad(\varOmega_0)+\mu(\varOmega')\right)^s
                                                            (\bnull \otimes \bunit). 
\end{align}
Here we use the following convention of notations:
\begin{align*}
\ad(\omega \otimes X)(\varphi\otimes F)&=(\omega\circ\varphi)\otimes \ad(X)(F),\\
\mu(\omega \otimes X)(\varphi\otimes F)&=(\omega\circ\varphi)\otimes \mu(X)(F)
\end{align*}
for $\varphi \otimes F \in S(A) \otimes \cU(\fX),\; \omega \otimes X \in \cB_1 \otimes \fX$.

This says that \textbf{the fundamental solution normalized at $(z_1,z_2)=(0,0)$ 
exists and is unique}. Moreover we can show that 
\begin{align*}\tag{IIFORM} \label{iiform}
    \Big(\ad(\varOmega_0)+\mu(\varOmega')\Big)^s(\bnull \otimes \bunit) \ 
    \in \cB^0 \otimes \cU_s(\fX).
\end{align*}

\section{Decomposition theorem and hyperlogarithms}
\subsection{The decomposition theorem of the normalized fundamental solution}

We consider the following four formal (generalized) 1KZ equation. 
In the following $d_{z_1}$ (resp. $d_{z_2}$) stands for the exterior differentiation 
by the variable $z_1$ (resp. $z_2$): 

\begin{alignat*}{4}
d_{z_1}G(z_1,z_2)&=\varOmega_{1\otimes2}^{(1)} G(z_1,z_2), \quad 
   &\varOmega_{1\otimes2}^{(1)}&=\zeta_1 Z_1+\zeta_{11} Z_{11}+\zeta_{12}^{(1)} Z_{12}, \\
d_{z_2}G(z_2)&=\varOmega_{1\otimes2}^{(2)} G(z_2), \quad 
   &\varOmega_{1\otimes2}^{(2)}&=\zeta_2 Z_2+\zeta_{22} Z_{22},\\
d_{z_2}G(z_1,z_2)&=\varOmega_{2\otimes1}^{(2)}G(z_1,z_2), \quad 
   &\varOmega_{2\otimes1}^{(2)}&=\zeta_2 Z_2+\zeta_{22} Z_{22}+\zeta_{12}^{(2)} Z_{12},\\
d_{z_1}G(z_1)&=\varOmega_{2\otimes1}^{(1)} G(z_1), \quad 
   &\varOmega_{2\otimes1}^{(1)}&=\zeta_1 Z_1+\zeta_{11} Z_{11}. 
\end{alignat*}
The fundamental solution normalized at the origin to each equation satisfies the conditions 
\begin{align*}
& \cL_{i_1 \otimes i_2}^{(i_k)} = \hcL_{i_1 \otimes i_2}^{(i_k)} \, z_{i_k}^{Z_{i_k}},\\
\hcL_{i_1 \otimes i_2}^{(i_k)} =\sum_{s=0}^{\infty} \hcL_{i_1 \otimes i_2,s}^{(i_k)},
\quad  & \hcL_{i_1\otimes i_2,s}^{(i_k)}\Big|_{z_{i_k}=0}=0 \quad (s > 0), \quad 
\hcL_{i_1\otimes i_2,0}^{(i_k)}= \bunit. 
\end{align*}

\begin{prop}\label{prop:decomposition}
\begin{enumerate}
\item The fundamental solution $\cL(z_1,z_2)$ of \eqref{2KZeq} normalized at the origin 
      decomposes to product of the normalized fundamental solutions of the (generalized) 
      formal 1KZ equations as follows: 
\begin{align*}
        \cL(z_1,z_2)&=\cL_{1 \otimes 2}^{(1)}\cL_{1 \otimes 2}^{(2)} 
                     = \hcL_{1 \otimes 2}^{(1)}\hcL_{1 \otimes2 }^{(2)}z_1^{Z_1}z_2^{Z_2} \\
                    &= \cL_{2 \otimes 1}^{(2)}\cL_{2 \otimes 1}^{(1)} 
                     = \hcL_{2 \otimes 1}^{(2)}\hcL_{2 \otimes 1}^{(1)}z_1^{Z_1}z_2^{Z_2}.
\end{align*}
\item If the decomposition 
\begin{align*}
       \cL(z_1,z_2)=G_{i_1\otimes i_2}^{(i_1)}G_{i_1\otimes i_2}^{(i_2)}
\end{align*}
      holds, where $G_{i_1\otimes i_2}^{(i_k)}=\hat{G}_{i_1\otimes i_2}^{(i_k)} \, z_{i_k}^{Z_{i_k}}$ 
      satisfies the same conditions as $\cL_{i_1 \otimes i_2}^{(i_k)}$ does, we have 
      $G_{i_1\otimes i_2}^{(i_k)}=\cL_{i_1 \otimes i_2}^{(i_k)}$.
\end{enumerate}
\end{prop}

\subsection{The iterated integral solution along the contours $C_{1 \otimes 2}$ and 
$C_{2 \otimes 1}$}

From \eqref{iiform}, we can choose $C_{1 \otimes 2}$ as the integral contour 
in \eqref{iisol}. Hence we have 
\begin{align*}
\hcL_s(z_1,z_2)&=\int_{C_{1 \otimes 2}} 
              \left(\ad(\varOmega_0)+\mu(\varOmega')\right)^s(\bnull \otimes \bunit)\\
           &=\int_{C_{1 \otimes 2}} (\iota_{1 \otimes 2}\otimes \id_{\cU(\fX)})
              \left(\left(\ad(\varOmega_0)+\mu(\varOmega')\right)^s(\bnull \otimes \bunit)\right) \\
           &=\sum_{s'+s''=s} \ \sum_{W',W''}
           \int_0^{z_1} \theta^{(1)}_{1 \otimes 2}(W')
              \int_0^{z_2} \ \theta^{(2)}_{1 \otimes 2}(W'') \ \alpha(W')\alpha(W'')(\bunit).
\end{align*}
Here $W'$ runs over $\cW_{s'}^0(Z_1,Z_{11},Z_{12})$, 
$W''$ runs over $\cW_{s''}^0(Z_2,Z_{22})$. \ 
($\cW^0_s(\fA)=\cW^0(\fA) \cap \cU_s(\fX)$, and $\cW^0(\fA)$
stands for the set of words of the letters $\fA$ which do not end with $Z_1,Z_2$.) 
$\alpha: \cU(\fX) \to \End(\cU(\fX))$ is an algebra homomorphism
\begin{align*}
    \alpha: (Z_1,Z_{11},Z_2,Z_{22},Z_{12}) \mapsto 
                   (\ad(Z_1),\mu(Z_{11}),\ad(Z_2),\mu(Z_{22}),\mu(Z_{12})),
\end{align*}
and $\theta^{(1)}_{1 \otimes 2}:\cU(\C\{Z_1,Z_{11},Z_{12}\}) \to S(\zeta_1,\zeta_{11},\zeta_{12}^{(1)})$ and $\theta^{(2)}_{1 \otimes 2}:\cU(\C\{Z_2,Z_{22}\}) \to S(\zeta_2,\zeta_{22})$ are linear maps
defined by replacing 
\begin{align*}
     \theta^{(i)}_{1 \otimes 2}(Z_i)=\zeta_i, \ \theta^{(i)}_{1 \otimes 2}(Z_{ii})=\zeta_{ii} \ \ (i=1,2),
      \ \theta^{(1)}_{1 \otimes 2}(Z_{12})=\zeta_{12}^{(1)}.
\end{align*}

\noindent
In the same way, we have 
\begin{align*}
  \hcL_s(z_1,z_2) &=\int_{C_{2 \otimes 1}} 
                   \left(\ad(\varOmega_0)+\mu(\varOmega')\right)^s(\bnull \otimes \bunit)\\
                  & =\int_{C_{2 \otimes 1}} (\iota_{2 \otimes 1}\otimes \id_{\cU(\fX)})
     \left(\left(\ad(\varOmega_0)+\mu(\varOmega')\right)^s(\bnull \otimes \bunit)\right)\\
     &=\sum_{s'+s''=s} \ \sum_{W',W''}
           \int_0^{z_2} \theta^{(2)}_{2 \otimes 1}(W')
              \int_0^{z_1} \ \theta^{(1)}_{2 \otimes 1}(W'') \ \alpha(W')\alpha(W'')(\bunit).
\end{align*}
Here $W'$ runs over  $\cW_{s'}^0(Z_2,Z_{22},Z_{12})$, and $W''$ runs over 
$\cW_{s''}^0(Z_1,Z_{11})$. 
$\theta^{(2)}_{2 \otimes 1}:\cU(\C\{(Z_2,Z_{22},Z_{12}\}) \to S(\zeta_2,\zeta_{22},\zeta_{12}^{(2)})$ and $\theta^{(1)}_{2 \otimes 1}:\cU(\C\{Z_1,Z_{11}\}) \to S(\zeta_1,\zeta_{11})$ are linear maps
defined by replacing 
\begin{align*}
     \theta^{(i)}_{2 \otimes 1}(Z_i)=\zeta_i, \ \theta^{(i)}_{2 \otimes 1}(Z_{ii})=\zeta_{ii} \ \ (i=1,2),
      \  \theta^{(2)}_{2 \otimes 1}(Z_{12})=\zeta_{12}^{(2)}.
\end{align*}

Since $[Z_1, Z_2]=[Z_1,Z_{22}]=0$, we have 
\begin{align*}
  \hcL(z_1,z_2)= \left( \sum_{W'}\int_0^{z_1} \theta^{(1)}_{1 \otimes 2}(W')\alpha(W') (\bunit) \right)
  \, \left( \sum_{W''}\int_0^{z_2} \theta^{(2)}_{1 \otimes 2}(W'')\alpha(W'') (\bunit) \right). 
\end{align*}
This says that \textbf{each decomposition 
in Proposition \ref{prop:decomposition} corresponds to 
the choice of the integral contours} $C_{1 \otimes 2}, \, C_{2 \otimes 1}$.

\subsection{Hyperlogarithms of the type $\cM_{0,5}$}

In \eqref{hlog}, let $m=2, a_1=1, a_2=z_2$, replace $\xi_0, \xi_1, \xi_2$ by $\zeta_1, \zeta_{11}, 
\zeta_{12}^{(1)}$ respectively, and put $\zeta(a_i)=\xi_i \ (i=1,2)$. 
Then \eqref{hlog} reads as
\begin{align*}
         L({}^{k_1} a_{i_1} \cdots{}^{k_r} a_{i_r} \, ;z_1) 
               &=\int_0^{z_1} \zeta_1^{k_1-1}\circ \zeta(a_{i_1}) \circ 
                                 \zeta_1^{k_2-1}\circ \zeta(a_{i_2})
                            \circ \cdots \circ \zeta_1^{k_r-1} \circ \zeta(a_{i_r}) \\
               & = \sum_{n_1>n_2>\cdots>n_r>0} \ \frac{a_{i_1}^{n_1-n_2} a_{i_2}^{n_2-n_3}
                            \cdots a_{i_r}^{n_r}}{n_1^{k_1}\cdots n_r^{k_r}} z_1^{n_1}, 
\end{align*}
which is referred to as a \textbf{hyperlogarithm of the type $\cM_{0,5}$}.
If $a_{i_1}=\cdots=a_{i_r}=1$, it is a multiple polylogarithm of one variable \eqref{1mpl}
\begin{align*}
        \Li_{k_1,\ldots,k_{r}}(z_1)=L({}^{k_1}1\cdots{}^{k_r}1;z_1), 
\end{align*}
and 
\begin{align}\tag{2MPL} \label{2mpl}
      \Li_{k_1,\ldots,k_{i+j}}(i,j;z_1,z_2)
                := L({}^{k_1}1\cdots{}^{k_i}1{}^{k_{i+1}}z_2 \cdots {}^{k_{i+j}}z_2;z_1)
\end{align}
is called a \textbf{multiple polylogarithm of two variables}. They constitute a subclass of 
hyperlogarithms of the type $\cM_{0,5}$. 

We should note that, in the previous subsection, the iterated integral  
\begin{align*}
     L\big( \theta^{(1)}_{1 \otimes 2}(W') \,;\,  z_1 \big) :=
     \int_0^{z_1} \theta^{(1)}_{1 \otimes 2}(W') \quad (W' \in \cW_{s'}^0(Z_1,Z_{11},Z_{12}))
\end{align*}
is a hyperlogarithm of the type $\cM_{0,5}$, 
and the iterated integral 
\begin{align*}
    L\big( \theta^{(2)}_{1 \otimes 2}(W'') \,;\, z_2 \big) := 
     \int_0^{z_2} \ \theta^{(2)}_{1 \otimes 2}(W'') \quad ( W''\in \cW_{s''}^0(Z_2,Z_{22}))
\end{align*}
is a multiple polylogarithm of one variable. Thus, \textbf{the normalized fundamental solution 
$\cL(z_1,z_2)$ is a generating function of hyperlogarithms of the type $\cM_{0,5}$}.

\section{Relations of multiple polylogarithms}
\subsection{Generalized harmonic product relations of hyperlogarithms}

From Proposition \ref{prop:b0}, one can define
\begin{align*}
\varphi(W',W'')=\iota_{1\otimes2}^{-1}(\theta_{1\otimes2}^{(1)}(W') 
                          \otimes\theta_{1\otimes2}^{(2)}(W'')) \in \cB^0
\end{align*}
for $W' \in \cW^0(Z_1,Z_{11},Z_{12}), \ W'' \in \cW^0(Z_2,Z_{22})$. 
Then we have
\begin{align*}
\int_{C_{1\otimes2}}\!\!\iota_{1\otimes2}(\varphi(W',W''))
=L(\theta_{1\otimes 2}^{(1)}(W');z_1)L(\theta_{1\otimes 2}^{(2)}(W'');z_2),
\end{align*}
and
\begin{align*}
\hcL_s(z_1,z_2)=\sum_{s'+s''=s}\ 
\sum_{\substack{W' \in \cW^0_{s'}(Z_1,Z_{11},Z_{12})\\
                                 W'' \in \cW^0_{s''}(Z_2,Z_{22})}} 
\int_{(0,0)}^{(z_1,z_2)} \ \varphi(W',W'') \ \alpha(W')\alpha(W'')(\bunit).
\end{align*}

\noindent
Since $\{\alpha(W')\alpha(W'')(\bunit) \,|\, W' \in \cW^0(Z_1,Z_{11},Z_{12}), \
W'' \in \cW^0(Z_2,Z_{22})\}$ is a linearly independent set, we obtain the following proposition:
\begin{prop}
We have 
\begin{align}\tag{GHPR} \label{ghpr}
     L(\theta^{(1)}_{1\otimes 2}(W');z_1)L(\theta^{(2)}_{1\otimes 2}(W'');z_2)
            =\int_{C_{2\otimes 1}} \ \iota_{2\otimes 1}(\varphi(W',W'')) 
\end{align}
for $W' \in \cW^0(Z_1,Z_{1 1},Z_{12}), W'' \in \cW^0(Z_2,Z_{22})$. 
\end{prop}
We call \eqref{ghpr} the \textbf{generalized harmonic product relations of hyperlogarithms}.

\begin{rem}
We have actually 
\begin{align*}
 \left(\ad(\varOmega_0)+\mu(\varOmega')\right)^s(\bnull \otimes \bunit)
                  =\sum_{s'+s''=s} \ \sum_{W',W''}  \ 
                               \varphi(W',W'') \otimes \alpha(W')\alpha(W'')(\bunit).
\end{align*}
For the proof, see \cite{OU1}. 
\end{rem}

\subsection{Harmonic product of multiple polylogarithms}

For $W'=Z_1^{k_1-1}Z_{11}\cdots Z_1^{k_i-1}Z_{11} Z_1^{k_{i+1}-1}Z_{12} 
\cdots Z_1^{k_{i+j}-1}Z_{12}, \ W''=\bunit$, we have 
\begin{align*}
         \int_{C_{1 \otimes 2}}\varphi(W',\bunit)=\Li_{k_1,\ldots,k_{i+j}}(i,j;z_1,z_2). 
\end{align*}
Hence \eqref{ghpr} for this case reads as 
\begin{align*}
        \Li_{k_1,\ldots,k_{i+j}}(i,j;z_1,z_2)=\int_{C_{2 \otimes 1}} \varphi(W',\bunit). 
\end{align*}
Moreover, by induction, one can prove that 
\textbf{the generalized harmonic product relations properly contain 
the harmonic product of multiple polylogarithms} such as \eqref{hpmpl}. 

Taking the limit, we have harmonic product of multiple zeta values. 
Thus \textbf{we can interpret the harmonic product of multiple zeta values 
as a connection problem for the formal KZ equation} such as \eqref{hpmzv}.

\section{The five term relation for the dilogarithm}

We define the action of $\fS_n$ on $\cM_{0,n}$ by $\sigma(x_i)=x_{\sigma(i)}$.
For $n=5$, the action of $\sigma=(23)(45) \in \fS_5 $ is given, in the cubic coordinates,  
by a birational transformation on $\bP^1\times\bP^1$ such as 
\begin{align*}
    \sigma(z_1,z_2) = \left(\frac{-z_1(1-z_2)}{1-z_1}, \, \frac{-z_2(1-z_1)}{1-z_2} \right). 
\end{align*}
It satisfies $\sigma^2=\id$ and preserves the divisors $D(\cM^{cubic}_{0,5})$.

Let $\sigma^* : \cB \to \cB$ be the pull back induced by $\sigma$, 
\begin{align*}
 \begin{array}{l}
 \sigma^*\zeta_1=\zeta_1+\zeta_{11}-\zeta_{22}, 
                            \quad \sigma^*\zeta_{11}=-\zeta_{11}+\zeta_{12}, \\
 \sigma^*\zeta_2=-\zeta_{11}+\zeta_2+\zeta_{22}, 
          \quad \sigma^*\zeta_{22}=-\zeta_{22}+\zeta_{12}, \quad \sigma^*\zeta_{12}=\zeta_{12}.
\end{array}
\end{align*}
and define an automorphism $\sigma_* : \cU(\fX) \to \cU(\fX)$ by
\begin{align*}
     (\sigma^* \otimes \id) \varOmega = (\id \otimes \sigma_*) \varOmega. 
\end{align*}
Hence we have 
\begin{align*}
\begin{array}{l}
     \sigma_*Z_1 = Z_1, \quad  \sigma_*Z_{11} = Z_1 - Z_{11} -Z_2, \\
     \sigma_*Z_2=Z_2, \quad \sigma_*Z_{22} = - Z_1 + Z_2 - Z_{22}, 
                  \quad \sigma_*Z_{12} = Z_{11} + Z_{22} + Z_{12}.
\end{array} \notag
\end{align*}

Since $(\id \otimes \sigma_*)^{-1} (\sigma^* \otimes \id) \varOmega = 
(\sigma^* \otimes \sigma_*^{-1}) \varOmega = \varOmega$, the function
\begin{align*}
\tcL(z,w)=(\sigma^* \otimes \sigma_*^{-1})\cL(z_1,z_2)=\cL(\sigma(z_1,z_2))
         \Big|_{Z \to \sigma_*^{-1}Z,\ (Z=Z_1,Z_{11},Z_2,Z_{22},Z_{12})}
\end{align*}
is also a fundamental solution of the KZ equation of two variables
which has the asymptotic behavior
\begin{align*}
     \tcL(z_1,z_2) \sim \bunit \left(\frac{-z_1(1-z_2)}{1-z_1}\right)^{Z_1} 
                    \left(\frac{-z_2(1-z_1)}{1-z_2}\right)^{Z_2} \qquad  (z_1,z_2) \to (0,0).
\end{align*}
Therefore the \textbf{connection formula} for $\cL(z_1,z_2)$ and $\tcL(z_1,z_2)$ is written as
\begin{align*}
        \tcL(z_1,z_2) &= \mathcal{L} (z_1,z_2) \exp(-\mathrm{sgn}(\mathrm{Im} z_1) \, \pi i Z_1) 
                         \, \exp(-\mathrm{sgn}(\mathrm{Im} z_2)\, \pi i Z_2).
\end{align*}
For the later use, it is convenient to rewrite this as follows: 
\begin{align*}
     (\sigma^*\cL)(z_1,z_2) &= (\sigma_*\cL)(z_1,z_2) 
             \exp(-\mathrm{sgn}(\mathrm{Im} z_1)\, \pi i Z_1) 
                  \, \exp(-\mathrm{sgn}(\mathrm{Im} z_2)\, \pi i Z_2).    
\end{align*}

\noindent
The terms $[Z_1,Z_{11}]$ and $[Z_2,Z_{22}]$ in the both sides above  
appear in $\sigma^*\hcL_2(z_1,z_2)$ and $\sigma_*\hcL_2(z_1,z_2)$. 
Comparing the coefficients of $[Z_1,Z_{11}]$, we have 
\begin{align}
             \Li_2\left(\frac{-z_1(1-z_2)}{1-z_1}\right)
             = \Li_{1,1}(1,1;z_1,z_2)-\Li_2(z_1)-\Li_{1,1}(z_1)+\Li_2(0,1;z_1,z_2), 
             \tag{L1} \label{landen1}
\end{align}
and comparing the coefficients of $[Z_2,Z_{22}]$, 
\begin{equation}
           \Li_2\left(\frac{-z_2(1-z_1)}{1-z_2}\right)
           = -\Li_{1,1}(1,1;z_1,z_2)-\Li_2(z_2)-\Li_{1,1}(z_2)+\Li_1(z_2)\Li_1(z_1).
           \tag{L2}\label{landen2}
\end{equation}
We should observe that \eqref{landen1} is regarded as a \textbf{``two-variables'' 
analogue of the Landen formula for  the dilogarithm} \cite{Le}. 
Since $\Li_2(0,1;z_1,z_2)=\Li_2(z_1z_2)$ and 
$\Li_{1,1}(z)=\frac{1}{2}\log^2(1-z)$, $\eqref{landen1} + \eqref{landen2}$ 
gives the \textbf{five term relation for the dilogarithm} \eqref{5term}:
\begin{align*}
  \Li_2(z_1z_2)= \Li_2 & \left(\frac{-z_1(1-z_2)}{1-z_1}\right)  + 
                \Li_2\left(\frac{-z_2(1-z_1)}{1-z_2}\right)  \notag \\
     & \qquad + \Li_2(z_1)+\Li_2(z_2)+\frac{1}{2}\log^2\left(\frac{1-z_1}{1-z_2}\right).
\end{align*}

\end{document}